\documentclass[12pt]{article}
\usepackage{amssymb}
\usepackage{amsmath,amsfonts,amsthm}
\usepackage{srcltx}
\usepackage{footnote}
\begin{document}                                 
	
	   \begin{center} \bf{ POWERS Vs. POWERS}\\ Pramod K. Sharma\\ e-mail: pksharma1944@yahoo.com\\
                                     	Department of Matematics, Sikkim University , Gangtok\\INDIA
                                 \end{center}
                                 
                                 \begin{center} \bf{ABSTRACT}  \end{center} 

                                            Let $ A \subset B$ be rings. An ideal $ J \subset B$ is called power stable in $A$ if $ J^n \cap A = (J\cap A)^n$ for all $ n\geq 1$. Further, $J$ is called ultimately power stable in $A$ if $ J^n \cap A = (J\cap A)^n$ for all $n$ large i.e., $ n \gg 0$. In this note, our focus is to study these concepts for pair of rings $ R \subset R[X]$ where $R$ is an integral domain. Some of the results we prove are: A maximal ideal $\textbf{m}$ in $R[X]$ is power stable in $R$ if and only if $ \wp^t $ is $ \wp-$primary for all $ t \geq 1$ for the prime ideal $\wp = \textbf{m}\cap R$. We use this to prove that for a Hilbert domain $R$, any radical ideal in $R[X]$ which is a finite intersection of G-ideals is power stable in $R$. Further, we prove that if $R$ is a Noetherian integral domain of dimension 1 then any radical ideal in $R[X] $ is power stable in $R$, and if every ideal in $R[X]$ is power stable in $R$ then $R$ is a field. We also show that if $ A \subset B$ are Noetherian rings, and  $ I $ is an ideal in $B$ which is ultimately power stable in $A$, then if $ I \cap A = J$ is a radical ideal generated by a regular $A$-sequence, it is power stable. Finally, we give a relationship in power stability and ultimate power stability using the concept of reduction of an ideal (Theorem 3.22). \\
                                                                             
                      \section{INTRODUCTION} 
											
				This article is largely based on forgotten note [6]. We assume all rings are commutative with identity. For a subset $S$ of a ring $R, id.(S)$ shall  denote the ideal of $R$ generated by $S$ and for an ideal $J$ of $R$, $\widehat{R}_J$ shall denote the $J$-adic completion of the ring $R$. If $I$ is an ideal of $R[X]$, then $ \overline{a(X)}$ denotes the image of $a(X) \in R[X]$ in the quotient ring $ R[X]/I$. We shall use $ \subset$ to mean contained or equal to.  In [2], an ideal $I$ in a ring $R$ is defined almost prime if for all $ a,b \in R, ab \in I- I^2$ either $a\in I \mbox{ or } b\in I$. While trying to prove that all ideals in $\mathbb{Z}[X]$ are almost prime we required that for any ideal $I \subset \mathbb{Z}[X], I\cap \mathbb{Z}= n \mathbb{Z}$ implies $ I^2 \cap \mathbb{Z}= n^2 \mathbb{Z}$. This, however, was not true. This property seems interesting in itself and is the basis of our definitions of power stable and ultimately power stable ideals in a polynomial ring $R[X]$ (Definition 2.1), and for pair of rings $ A \subset B$ (Definition 2.2). In this article $R$ shall always denote an integral 
             domain. We prove that a maximal ideal $\textbf{m}$ in $R[X]$ is power stable if and 
             only if $ \wp^t $ is $ \wp$-primary for all $ t \geq 1$ for the prime ideal $\wp = \textbf{m}\cap R$ (Theorem 3.10). This result is used to prove that if $R$ is a Hilbert domain then any radical ideal in $R[X]$ which is a finite intersection of $G-$ideals is power stable (Theorem 3.14). Further, it is proved that if $R$ is a Noetherian domain of dimension 1, then any radical ideal in $R[X]$ is power stable (Theorem 3.15) and if every ideal in $R[X]$ is power stable then $R$ is a field (Theorem 3.18) . We also prove that if $ A \subset B$ are Noetherian rings, and  $ I $ is an ideal in $B$ which is ultimately power stable in $A$, then if $ I \cap A = J$ is a radical ideal generated by a regular $A-$ sequence, it is power stable. Finally, we give a relationship in power stability and ultimate power stability using the concept of reduction of an ideal (Theorem 3.22). \\

              \section{Observations and Definitions}       
             
             Definition 2.1. Let $R$ be an integral domain and $I \subset R[X]$, an ideal in the polynomial ring $R[X]$. Then\\
             
             (i) The ideal $I$  is called power stable if for all $t\geq 1, I^t\cap R = (I\cap R)^t$.\\  
                        
             (ii) The ideal $I$ is called ultimately power stable if $ I^t\cap R = (I\cap R)^t$ for all $ t \gg 0$. \\
             
             More generally, we define:\\ 
             
             Definition 2.2.  If $ A \subset B$ are rings, then an ideal $J$ in $B$ is called power stable (respectively  ultimately power stable) in $A$ if $ J^n \cap A = (J\cap A)^n$ for all $ n \geq 1$ (for all $ n \gg 0$). \\
             
             Let us first of all note that an ultimately power stable ideal need not be power stable in general. Consider $A = \mathbb{Z}/(4) \subset B = \mathbb{Z}/(4)[X]/(X^2- 2)$. Then the ideal $I = (\overline{2}, \overline{X} ) $ in the ring $B$ satisfies $ I\cap A = I^2\cap A = (\overline{2})$. Thus  $I^2 \cap A \neq (I\cap A)^2$. However $ I^m \cap A = (0)= (I\cap A)^m$ for all $ m \geq 3$. Thus the ideal $I$ in the ring $B$ is not power stable in $A$, but is ultimately power stable. We would like to know if the same holds in $R[X]$, but at present we have no example to this effect. \\
                                        
             Example 2.3. Any principal ideal in $R[X]$ is power stable.\\
             
             Example 2.4. Let $ A \subset B$ be rings. Then any ideal $ J \subset B$ with $ J \cap A =(0)$ is power stable in $A$.\\ 
             
             Example 2.5. For any ideal $I \subset R, I[X]$ is power stable.\\
             
             Example 2.6.  Let $ A \subset B$ be rings. If an ideal $J \subset B$ is power stable in $A$, then $J^t$  power stable ideal in $A$ for all $ t \geq 1$. \\
             
             Example 2.7. Let $ A \subset B$ be rings. If $J \subset B $ is an ideal power stable  in $A$ and $ J\cap A = I$, then for any ideal $J_1 \subset J$ in $B$ with $ J_1 \cap A = I$, $ J_1$ is power stable in $A$.\\
             
             Example 2.8. Let $ A\subset B$ be rings where  $B$ is faithfully flat over $A$. Then for any ideal $I \subset A, IB$ is power stable in $A$. .\\
             
             Example 2.9. Let $ A\subset B$ be rings. If $B$ is flat as well as integral over $A$, then for any ideal $ I \subset A, IB$ is power stable over $A$. \\
             
             First of all, we make a few general observations.\\
             
             Lemma 2.10. If $\varphi$ is an automorphism of $ R[X]$ such that $ \varphi (R) = R$, then for any power stable ideal $I \subset R[X], \varphi(I)$ is power stable.\\  
             
             Proof. It is clear.\\         
                            
             Lemma 2.11. Let $ A \subset B$ be rings, and  $I$ be an ideal  in $B$ where $ I\cap A = J$. Then \\ (i) The ideal $I$ is power stable in $A$ if and only if the natural homomorphism $$ \varphi : Gr_J (A) \longrightarrow  Gr_I(B) $$ is a monomophism of graded rings.\\
             
             (ii) If $I$ is ultimately power stable in $A$, and $ Gr_J(A)$ has no nilpotents of degree $ \geqq 1$, then $I$ is power stable in $A$.\\
                          
             Proof. (i) If $I$ is power stable in $A$, then $ I^n \cap A = J^n$ for all $ n \geqq 0.$ Hence
             
       \begin{center}
       	      \begin{eqnarray*} J^n\cap I^{n+1} & = & I^n\cap A \cap I^{n+1}\\
             	 & = & I^{n+1} \cap A\\
             	 & = & J^{n+1}            
             \end{eqnarray*}      
  
       \end{center}           
             Therefore $\varphi$ is a monomorphism. Conversly, let $\varphi$ be a monomorphism. Then \\ \begin{center} $ J^n \cap I^{n+1} = J^{n+1} $\hspace*{2 cm} (1)   \end{center}
              for all $ n \geqq 0$. We shall prove that $ I^n \cap A = J^n $ for all $ n \geqq 1$ by induction on n. Since $ \varphi $ is a monomorphism, the statement is clear for $ n= 1.$
              Let $ n \geqq 2$. By induction assumption, 
              
                 \begin{eqnarray*}
                 	 I^{n-1}\cap A  & = &J^{n-1}\\ 
      \Rightarrow I^n\cap (I^{n-1} \cap A) & = &  I^n \cap J^{n-1}\\ 
      \Rightarrow I^n \cap A & = & J^n                         
                 \end{eqnarray*}
      by equation (1). Hence the result follows.     \\
      
      (ii) Note that  the natural morphism $\varphi$ from $ Gr_J(A)$ to  $ Gr_I(B)$ is a homomorphism of  graded rings and has degree zero. As $I$ is ultimately power stable $\varphi$ is monomorphism on homogeneous components of degree $n \gg 0$. Clearly as $Gr_J(A)$ has no nilpotents of degree $ \geqq 1, \varphi$ is injective. Thus by (i), the ideal $I$ is power stable in $A$.\\ 
      
      Lemma 2.12. Let $ A \subset B$ be rings, and $I$ be an ideal of $ B$. Let  $ I\cap A = J$. Then \\ (i) If $I$ is power stable in $A$, the natural homomorphism  $$ \widehat{A}_J \longrightarrow  \widehat{B}_I$$  is a
       monomorphism. \\ (ii) If the rings $A, B$ are Noetherian, then $I$ is power stable in $A$ if and only if  the natural homomorphism 
      $$ \widehat{A}_J/ (J^n) \longrightarrow  \widehat{B}_I/ (I^n)$$
      is a  monomorphism for all $ n\geqq 1$.\\
      
      Proof. First of all, note that $ \widehat{A}_J = \varprojlim A/J^n $ and $\widehat{B}_I = \varprojlim B/I^n.$ We, now, prove :\\ (i) As $I$ is power stable, $I^n\cap A = J^n$ for all $ n \geq 1$. Hence the natural map  $ A/J^n \stackrel{\alpha_n}\longrightarrow B/I^n $ is a monomorphism for all n, and the diagram :\\
      
                               $$  \begin{array}{ccc} A/J^{n+1}  &   \stackrel{\alpha_{n+1}} \longrightarrow & B/I^{n+1} \\ \downarrow &  & \downarrow \\ A/J^n & \stackrel{\alpha_{n}} \longrightarrow & B/I^n  \end{array} $$ is commutative where the vertical maps are quotient maps. This set up clearly induces a natural monomorphism $ \widehat{A}_J \longrightarrow \widehat{B}_I.$\\ (ii) In case $A,B$ are Noetherian, the diagram : $$  \begin{array}{ccc} \widehat{A}_J/ (J^n)  &   \stackrel{\beta_{n}} \longrightarrow & \widehat{B}_I/I^{n} \\ \uparrow &  & \uparrow \\ A/J^n & \stackrel{\alpha_{n}} \longrightarrow & B/I^n  \end{array} $$ is commutative where all the morphisms are natural. The vertical  maps are isomorphisms. Thus it is clear that $\alpha_n$ is a monomorphism if and only if $\beta_n$ is a monomorphism.\\
                               
                               Lemma 2.13. Let $ A \subset B\subset C $ be rings. If an ideal $J \subset C$ is power stable (ultimately power stable) in $B$ and $ J\cap B$ is power stable ( ultimately power stable)in $A$, then $J$ is power stable ( ultimately power stable) in $A$.\\
                               
                               Proof. It is straight forward.\\ 
                                 
                               Lemma 2.14. Let $ A \subset B$ be rings, and let $ I, J$ be ideals in $B$. If $ I,J$ are power stable in $A$ and $ I\cap A, J \cap A$ are comaximal, then $IJ = I\cap J$ is power stable in $A$.\\
                               
                               Proof. First of all, it is clear that $$ (IJ\cap A)^n \subset (IJ)^n \cap A$$ for all $ n \geq 1$. Now as $ I\cap A \mbox { and  } J\cap A$ are comaximal, $ I, J$ are comaximal in $B$. Hence
                                \begin{eqnarray*} (IJ)^n \cap A & = & (I^n \cap A )\cap (J^n \cap A)  \\
                               	\Longrightarrow   (IJ)^n \cap A & = & (I\cap A)^n \cap (J\cap A)^n \mbox {	since I and J are power stable in A.} \\
                               	  & =  &  ( (I\cap A )(J \cap A) )^n  \mbox { since } (I\cap A) \mbox{ and } (J\cap A) \mbox{ are co-maximal }.\\
                               	   &  \subset   &  (IJ\cap A)^n.                               	   
                               	  \end{eqnarray*}
                                 Consequently $IJ$ is power stable in $A$.  \\

                                Theorem 2.15.  Let $I$ be an ideal in a ring $A$, then there exists a ring $ A \subsetneqq B$ and a principal ideal $J$ in $B$ such that $J$ is power stable in $A$ and $J \cap A =I$. \\
                                
                                Proof. The generalised Rees ring $ B = A[ It, t^{-1}]$ does the job taking $ J = id.(t^{-1})$.\\
                                
                                Theorem 2.16. Let $ A \subsetneqq B$ be rings such that $B$ is flat over $A$. Then every ideal $I$ of $A$ is a contraction a power stable ideal of $B$ in $A$ if and only if $B$ is faithfully flat over $A$.\\
                                
                                Proof. If $I$ is contraction of a power stable ideal of $B$ in $A$, then $IB$, the extention of $I$ in $B$, is power stable with contraction $I$. Thus $B/IB$ is not zero, and hence for any finitely generated non-zero $A-$module $M$, $B\otimes_A M$ is not zero. Consequently $B$ is faithfully flat over $A$. The converse is already noted in the example 2.8.\\
                                
                                Theorem 2.17. Let $R$ be Noetherian. If every irreducible ideal in $R[X]$ is a finite intersection of power stable ideals, then every ideal is a finite intersection of power stable ideals.\\
                                
                                Proof. Let $$\mathcal{S}= \{ I \subset R[X],  \mbox { an ideal } : I \mbox { is not a finite intersection of power stable ideals} \}.$$ If $ \mathcal{S}= \phi $, we have nothing to prove. If not, then let $ J $ be a maximal element in $ \mathcal{S}$. By assumption, $ J $ is not irreducible. Hence $ J = K \cap L$ where $ K, L$ are ideals in $ R[X]$ strictly bigger than $J$. Hence  $K \mbox { and } L $ are finite intersections of power stable iodeals, and cosequently $ A$ is a finite intersection of power stable ideals. This cotradicts the fact that $ A \in \mathcal{S}$. Hence the assertion follows.\\
                                
                                Theorem 2.18.  Let $ A \subset B$ be rings where $B$ is Noetherian. Let $I$ be a regular ideal in $B$ which is power stable in $A$. Then the ideal $ I^* = \cup \{ I^{(n+1)} : I^n ; n\geqq 1 \}$is ultimately power stable.\\ 
                                
                                Proof.

                               \section{ Main Results} 
                                 
                       Lemma 3.1. An ideal $I$ in $R[X]$ is power stable if and only if $I_{\wp}$ is power stable in $R_{\wp}[X] \text { for all } \wp \in Spec(R)$.\\
                      
                      Proof. If $I$ is power stable then $ I_{\wp}$ is power stable for any prime ideal $ \wp  \in Spec(R)$ since localization commutes with intersections and powers. Note that $ ( I\cap R)^n \subset I^n\cap R \text{ for all } n \geqq 1$. Further, if $I_{\wp} $ is power stable for all primes  $\wp \in Spec(R), (I_{\wp})^n \cap R_{\wp} =  (I_{\wp}\cap R_{\wp})^n$. Consequently $ ((I^n\cap R)/( I\cap R)^n)_{\wp} = 0$ for all prime ideals  $\wp \in Spec(R)$. This implies $I^n\cap R \subset (I\cap R)^n$. Hence $ I^n\cap R = (I\cap R)^n$ for all $n \geqq 1$ i.e., $I$ is power stable.\\
                      
                      Remark 3.2.(i) As in above Lemma, if $ A\subset B$ are rings and $ I \subset B$ is an ideal, then $I$ is power stable in $A$  if and only if $I_{\wp}$  is power stable in $A_{\wp}$ for every prime ideal $\wp \in  Spec(A)$. Thus clearly it is sufficient to check power stability at maximal ideals in $A$.\\
                      
                      (ii) Let $R$ be an almost Dedekind domain. If $ I_1 \subsetneqq I$ be power stable ideals in $R[X]$ such that $ I_1 \cap R \subsetneqq I \cap R$ , then $ {I_1}^n \subsetneqq I^n$ for all $n \geqslant 1$.\\

                      Theorem 3.3. Let $R$ be a principal ideal domain, and for an ideal $I$ in $R[X], I\cap R = Rd$. If the image of $I$ in $R/(d)[X]$ is generated by a regular element then $I$ is power stable.\\
                      
                      Proof. If $ d = 0$, then the result is clear. Assume $ d\neq 0$. By assumption on $I, I = id.(d,h(X))$ for an $h(X) \in R[X]$, where image of $h(X)$ in $ R/(d)[X]$ is regular.We shall prove by induction  that $ I^t \cap R = (I \cap R)^t$  for all $ t \geqq 1$. Let $ t \geqq 2$, and $ I^s\cap R = (I\cap R)^s$ for all $ s \leqq t-1$. If $ I^t\cap R = eR$, then \\
                      
                      \begin{eqnarray*} eR & = & I^t\cap R  \subset  I^{t-1} \cap R \\ & = & d^{t-1}R \\ \Rightarrow e & = & d^{t-1}k \end{eqnarray*}
                       Now, as $ e\in I^t$, we have 
                      \begin{eqnarray*} d^{t-1}k & = &   \sum_{i+j=t, i\geqq 1} d^i h^j(X) a_{ij}(X) + h^t(X)a(X) \text{ for some } a_{ij}(X), a(X) \in R[X]. \\ \Rightarrow \overline{h^t(X) a(X)} & = & 0 \text{ in } R/(d)[X]\\ 
                      	\Rightarrow  \overline{a(X)} & = & 0 \text{ in } R/(d)[X] \text{ since } \overline{h(X)} \text{ is regular in } R/(d)[X]. \\ \Rightarrow a(X) & = & d b(X)  \text{ where }   b(X) \in R[X]\\ \Rightarrow d^{t-1} k & = & \sum d^i h^j(X)a_{ij}(X) + d h^t(X) b(X) \\ \Rightarrow d^{t-2} k & = & \sum d^{i-1} h^j(X) a_{ij}(X) + h^t(X) b(X) \\ \Rightarrow d^{t-2} k & \in & I^{t-1} \cap R \\   & = & d^{t-1}R \\ \Rightarrow k & \in & dR \\  \Rightarrow  e & = & d^t k_1 \text{ for some }  k_1 \in R                       	
                      	\end{eqnarray*}  Hence $I^t\cap R = d^t R$, and the result follows.\\
                      	
                      	Corollary 3.4. Let $R$ be a principal ideal domain. Then any prime ideal $\wp$ in $R[X]$ is power stable.\\
                      	
                      	Proof. As any non-zero prime ideal in $R$ is maximal, we have either $ \wp \cap R = (0) $ or a maximal ideal. Hence the proof is immediate from the theorem.\\ 
                      	
                      	Remark 3.5. In the Theorem 3.3 we can replace $R[X]$ by any integral domain $A$.\\

                      	Theorem 3.6. Let $R$ be a principal ideal domain, and  $ I \subset R[X]$, an ideal. Let  $ I = id.(d, h(X)) \text { such that }  d(\neq 0) \in R \text { and } h(X) \in R[X]$ is a non constant polynomial. Then $I $ is power stable.\\
                      	
                      	Proof. If $ h(X) \in R[X]d$, the result is clear. Hence, assume $ h(X) \notin R[X]d$.  Let $e$ be the greatest common divisor of $d$ and coefficients of $h(X)$ in R. Then $ I = e  (id.(d_1, h_1(X)) \text{ where } d= e d_1 \text { and } h(X) = e h_1(X)$. Put $ I_1 = id.(d_1, h_1(X)$) It is easy to check that for any ideal $J \subset R[X] \text { and } \lambda \in R, $ $ (\lambda J) \cap R = \lambda (J \cap R) $ . Hence, if we prove that $I_1$ is power stable then 
                      	
                      	\begin{eqnarray*}                      	
                      	  I^n \cap R & =  & e^n (I_1)^n \cap R \\
                      	  &  = & e^n ((I_1)^n\cap R ) \\ &  = & e^n (I_1 \cap R)^n \\ &  = &  ( e (I_1 \cap R) )^n \\ & = & (I \cap R )^n  \end{eqnarray*} i.e., $I$ is power stable.  We shall, now, show that $I_1$ is power stable. Note that $ \overline{h_1(X)}$ is a regular element in $R/(d_1)[X]$ since otherwise there exists $ \lambda \in R$ with its image $ \overline{\lambda}$ in $R/(d_1)$  non zero such that
                      	\begin{eqnarray*} 
                      		      \overline{\lambda. h_1(X)} & =  & 0  \text{ in } R/(d_1)[X] \\ \Rightarrow  \lambda h_1(X) &  \in & d_1R[X] \end{eqnarray*}
                   Hence, as greatest common divisor of $d_1$ and elements of $h_1(X)$ is identity,  $d_1 \text { divides }\lambda$, i.e., $ \overline{\lambda}= 0 \text { in }  R/(d_1)$, a contradiction to our assumption. Thus $\overline{h_1(X)}$ is regular and the result follows from Theorem 3.3. \\  
                   
                   Theorem 3.7. Let $ I \subset R[X]$ be an ideal such that $ I = id.( J, f(X))$ where $J$ is an ideal in $R$ and $ f(X) \in R[X]$ is a monic polynomial of degree $ \geqq 1$. Then $I$ is power stable.\\
                   
                   Proof. Let $\lambda \in I \cap R$, then $$  \lambda = a(X) + f(X) h(X)$$ where $ a(X) \in J[X]$. Reading off this equation in $ R/J[X]$, we get  \begin{eqnarray*} \overline{\lambda} & = & \overline{f(X) h(X)}\\ \Rightarrow \overline{\lambda} & = & 0 \text{ since } \overline{f(X)} \text { is monic of degree } \geqq 1 \text{ in } R/J[X]\\ \Rightarrow \lambda & \in & J \end{eqnarray*}
                   Hence $ I\cap R = J$. Now, let $ t > 1 \text{ and } \lambda \in I^t \cap R $, then $$ \lambda = a(X) + f(X) h(X) + f^t(X) c(X)$$ where $ a(X) \in J^t[X] \text{ and } h(X) \in J[X]$. As in case $t=1$, reading off this equation in $ R/J[X]$, we conclude $ c(X) \in J[X] \text{ and } \lambda \in J$. Hence $$ \lambda = a_1(X) + f(X) b_1(X)$$ where $a_1(X) \in J^t[X] \text{ and } b_1(X) \in J[X]$ where either $b_1(X) = 0$ or no coeficient of $b_1(X)$ lie in $J^t$. If $b_1(X) = 0$, then clearly $ \lambda \in J^t$. However, if $b_1(X) \neq  0 $, then if $ \lambda_0$ is leading coeficient of $a_1(X)$ and $ \mu_0$ is  is leading coeficient of $b_1(X) , \lambda_0 + \mu_0 = 0$. This implies $ \mu_0 \in J^t$. A contradiction to our assumption. Hence $ \lambda \in J^t$. therefore $I^t \cap R = J^t = (I \cap R)^t$ for all $ t \geqq 1$. Thus $I$ is power stable.\\
                   
                   Corollary 3.8. If for an ideal $ I \subset R[X], I\cap R = \textbf{m}$ is a maximal ideal, then $I$ is power stable.\\
                   
                   Proof. If $ I = \textbf{m}[X]$, then clearly $I$ is power stable. However, if $ I \neq \textbf{m}[X]$ then $ I = id.(\textbf{m}, f(X))$ where  $ f(X) \in R[X]$ is a monic polynomial of degree $ \geqq 1$. Therefore the result follows from the Theorem.\\  
                   
                   Corollary 3.9. Let $R$ be an integral domain of dimension one. Then any prime ideal of $R[X]$ is power stable.\\
                   
                   Proof. It is clear using Corollary 3.8.\\
                   
                   Theorem 3.10. A maximal ideal $ \textbf{m} \subset R[X]$ is power stable if and only if  for the prime ideal $ \wp =\textbf{m}\cap R, \wp^{(t)} = \wp^t \text { for all } t \geqq 1$, i.e., $\wp^t $ is $\wp-$primary  for all $ t \geqq 1$.\\                 
                     
                    Proof.  Let $ \textbf{m}$ be power stable. As $ \textbf{m}^t$ is $\textbf{m}-$primary $\text{ for all } t \geqq 1, \textbf{m}^t \cap R = \wp^t \text{ is }  \textbf{m}\cap R = \wp- $primary for all $ t\geqq 1$. Conversly, let $ \wp^t$ is $\wp-$primay for all $ t \geqq 1$. If $ \wp = (0)$, there nothing to prove. Hence let $ \wp \neq (0)$.Then, as $\textbf{m}_{\wp} \cap R_{\wp} =  \wp R_{\wp}, \textbf{m}_{\wp} \text{ is power stable by Corollary 3.8 }$. Hence
                      \begin{eqnarray*} 
                         (\textbf{m}_{\wp})^t \cap R_{\wp}& =  & \wp^t R_{\wp}\\ \Rightarrow \textbf{m}^t \cap R & \subset & \wp^{(t)} \\ \Rightarrow \textbf{m}^t \cap R & \subset & \wp^t \text{ since } \wp^{(t)} = \wp^t.\\ \Rightarrow   \textbf{m}^t \cap R  & = & (\textbf{m} \cap R)^t \text{ since }  \wp =\textbf{m}\cap R \text{ and }  \wp^t \subset \textbf{m}^t \cap R. \end{eqnarray*}  Thus the result is proved.\\ 	 
                                              
                    Remark 3.11. (1) In the reverse part of the above result, it is not used that $ \textbf{m}$ is maximal. Thus if $\textbf{q}$ is a prime ideal in $R[X]$ and for $ \textbf{q}_1 = \textbf{q} \cap R, \textbf{q}_1^t = \textbf{q}_1^{(t)}$ for all $ t \geq 1$, then $  \textbf{q}$ is power stable. Further, note that if $ \wp$ is a power stable prime ideal in $R[X]$, then for $\wp_1 = \wp \cap R, \wp_1^t $ need not be $\wp_1-$primary for all $ t \geq 1$. This follows since $\wp[X]$ is power stable in $R[X]$ for any $ \wp \in Spec(R)$. Thus if $ \wp^t$ is not  $\wp-$    primary for all $t \geq 1$, we get the required example.\\
                    
                    (2) If $R$ is a Hilbert domain then any maximal ideal in $R[X]$ is power stable. In particular, if $K$ is a field then for $ R=  K[X_1, \ldots, X_n]$, any maximal ideal in $R[X]$ is power stable.\\
                    
                    We, now, give two examples to show that, in general, a maximal ideal in $ R[X]$ need not be power stable. In view of Theorem 3.10, it suffices to give a G-ideal $P$ in $R$ for which $P^n$ is not $P-$primary for some $ n\geq 1$. The first example below was suggested by Melvin Hochster.\\
                    
                    Examples 3.12. Let $K$ be a field and $ Y,Z,W$ be algebraically independent over $K$. For an algebraically independent element $T$ over $K$, consider the $K-$algebra homomorphism: $$ \varphi: K[[Y,Z,W]] \longrightarrow K[[T]]$$ such that $\varphi(Y) =T^3 , \varphi(Z) = T^4, \varphi(W) = T^5$. Then kernel of $ \varphi$ is the prime ideal $ \wp = id.(f,g,h)  \text{ where }  f=(Y^3 - ZW), g = ( Z^2 - YW),\text { and } h = (W^2 - Y^2 Z) ) $. It is easy to see that $ \wp$ is a G-ideal. Further,we have $ f^2 - gh = Y p(Y,Z,W)\in \wp^2 $. Clearly $ Y \notin \wp $ and it is easy to check that  $ p(Y,Z,W) \notin \wp^2$. Thus $\wp^2$ is not $\wp -$primary. Hence for the integral domain $ R = K[[Y,Z,W ]]$, there is a maximal ideal $ \textbf{m}$ in the polynomial ring $ R[X]$ such that $ \textbf{m}\cap R = \wp$  and this $\textbf{m}$ is not power stable.\\
                    
                   Example 3.13. Let $K$ be a field, and $R = K[X,Z] [[Y]]/ (XY - Z^2)$. Note that $ K[X,Z][[Y]]$ is a unique factorization domain. Further, let   $$ Z^2-XY = (\sum_{ i\geqq 0} g_i(X,Z) Y^i ) (\sum_{ i\geqq 0} h_i(X,Z) Y^i ) $$ in $ K[X,Z][[Y]]$. Then, putting $ Y = 0$ in this equation, we get  $ Z^2 = g_0(X,Z). h_0(X,Z)$. Therefore $ g_0(X,Z) = \pm 1 \mbox{ or } \pm Z \mbox { or }  \pm Z^2$. Clearly if $ g_0(X,Z) = \pm 1 \mbox{ or } \pm Z^2$, then one of the factors in the equation is a unit. Thus let $ g_0(X,Z) = \pm Z$.  Consequently we also have  $ h_0(X,Z) = \pm Z$. Now, putting $ Z = 0$ in the equation, we get  $$ XY =   (\sum_{ i\geqq 1} g_i(X,0) Y^i)  (\sum_{ i\geqq 1} h_i(X,0) Y^i) $$ This, however is not possible. Hence $ Z^2 - XY$ is a prime element in $K[X,Z][[Y]]$. Thus $R$ is an integral domain. Now , let $ \wp = id.(\overline{X}, \overline{Z}) $ . Then $ \wp $ is a G-ideal in $R$, since $ R/\wp = k[[Y]]$ is a G-domain. Now note that $ \wp^2 = id.(\overline{X}^2, \overline{XZ}, \overline{Z}^2 )$ We have $ \overline{X}.\overline{Y}= \overline{Z}^2 \in \wp^2, \mbox{ but } \overline{X} \notin \wp^2, \mbox{ and } \overline{Y}\notin \wp $. Hence $\wp^2$ is not $\wp-$ primary. Thus there exists a maximal ideal $ \textbf{m}$ in the polynomial ring $ R[X]$ such that $ \textbf{m}\cap R = \wp$  and this $\textbf{m}$ is not power stable.\\

                    Theorem 3.14. Let $R$ be a Hilbert ring. Then any radical ideal $I \subset R[X]$, which is a finite intersection of G-ideals, is power stable.\\
                    
                    Proof. By [3, Theorem 31],  $R[X]$ is a Hilbert ring. Hence every G-ideal in $R[X]$ is maximal. Now, by our assumption on $I$, $$ I = \textbf{M}_1 \cap \textbf{M}_2 \ldots, \textbf{M}_n$$ where $ \textbf{M}_i$'s are distinct G-ideals in $R[X]$. As $R$ is a Hilbert ring $ \textbf{M}_i \cap R = \textbf{m}_i$ is a maximal ideal in $R$ for all $ i \geq 1$. Now note that for any $ t \geq 1$, \begin{eqnarray*} I^t & = & \textbf{M}^t_1 \cap \textbf{M}^t_2 \ldots, \textbf{M}^t_n \\ \Rightarrow I^t \cap R & = &  \textbf{M}^t_1 \cap \textbf{M}^t_2 \ldots, \textbf{M}^t_n \cap R \\ & = &  \textbf{m}^t_1 \cap \textbf{m}^t_2 \ldots, \textbf{m}^t_n       \end{eqnarray*} since, by Corollary 3.8, every maximal ideal in $ R[X]$ is power stable. Therefore it is clear that  $ I^t \cap R =( I \cap R)^t$ for every $ t\geq 1$ i.e., $I$ is power stable.\\ 
                    
                    Theorem 3.15. Let $R$ be a Noetherian domain of dimention 1.Then any radical ideal in $R[X]$ is power stable.\\
                    
                    Proof. If $I\cap R = (0)$, the result is clear. Hence, assume $ I \cap R = J \neq (0)$. Since $I$ is radical ideal of $R[X], J$ is a radical ideal of $R$. As $R$ is Noetherian of dimension 1, we have $$ J =  \textbf{M}_1 \cap \textbf{M}_2 \ldots, \textbf{M}_n$$ where $ \textbf{M}_{i}'s$ are maximal ideals in $R$.Thus it is clear that for any prime ideal $ \wp$ in $R$, either $ J_{\wp} = R_{\wp} \text{ or }  J_{\wp} = \wp_{\wp} $. Therefore, since $I_{\wp} \cap R_{\wp} = J_{\wp}$ for every prime ideal $ \wp$ in $R$, by Corollary 3.8, $I_{\wp}$ is power stable for all prime ideals $\wp$ in $R$. Hence by Lemma 3.1, $I$ is power stable. \\
                    
                    We shall now show that for an integral domain $R$ of of dimension 1, a non-radical ideal in $R[X]$ need not be power stable. In fact, we shall give an example of a primary ideal in $R[X]$, where $R$ is a principal ideal domain, which is not power stable. This  generalises an example given by Melvin Hochster who proved in a personal communication to Stephen McAdam that for any prime $p$, the ideal $ (X^2- p, X^3)$ in  $ \mathbb{Z}[X] $ is not power stable. We learned this  from Stephen McAdam. We note that this is not even ultimately power stable.\\
                    
                    Example 3.16. Let $p$ be a prime in a principal ideal domain $R$. Then for any two postive integers $m,n, n < m < 2n \mbox{ with } 2(m-n) \leq n, $ the ideal $ I =  (X^n- p, X^{m})$  is not power stable in $R[X]$.\\
                    
                    Proof. Step 1.  $I \cap R = Rp^2$.\\
                    
                                    We have $$ X^{m} - X^{m-n}(X^n-p)  =  pX^{m-n} \in I.$$  As $m-n < n, pX^n \in I .$ Hence $$ p(X^n - p) - pX^n = p^2 \in I$$
                                 	                                   
                                      If $ I \cap R = Rd$, then $ d$ divides $p^2$. We, now, note: \\
                                    
                                    (i) $ d\neq 1$\\
                                        If $ d=1$, then $I = R$. This, however, is not true as $ I \subset (X,p) $, and $(X,p) \cap R = Rp$. Thus $ d  \neq 1.$\\
                                        	 
                                 (ii)  $d \neq p.$\\
                                 
                                 Note that $ p\in I$ if and only if $ X^n \in I$. Thus if $d = p$, then $ X^n \in I$. Therefore we can write $$ X^n = (X^n - p) a(X) + X^{m} b(X) \hspace{1 cm }  ( a(X), b(X) \in R[X])$$  
                                 
                                 Substituting  $ X = 0$ in the above equation, we get $  a(0) = 0.$\\ Thus $a(X)$ is a multiple of $X$. Let $ a(X) = X^t a_1(X)$ such that $a_1(X)$ is not a multiple of $X$. Then  $$ X^n = (X^n - p)X^t a_1(X) + X^{m} b(X)$$ If $ t \gneqq n$ , then $$ 1 = (X^n - p)X^{t-n} a_1(X) +b(X) X^{m-n}$$ Putting $ X = 0$ in this equation gives $ 1 = 0$, an absurdity. Hence $ t \leq n$. In this case we get $$ X^{n-t} = a_1(X) (X^n - p) + b(X) X^{m-t}$$ Now, putting $X= 0$, we get $$ a_1(0) p = 0 \mbox{ or } -1.$$ None of this is possible since $p$ is a prime and $a_1(X)$ is not a multiple of $X$. Hence $p \notin I$. Thus step 1 is proved i.e., $I\cap R = Rp^2$ . \\ 
                                   
                                 Step 2. $ p^3 \in I^2 \cap R$.\\
                                 
                               Note that $$ I^2 = id.(X^{2n} -2p X^n + p^2, X^{2m}, X^{m+n} - p X^{m}).$$  Therefore as $ n < m$,  $$   X^{2m} - X^{m-n}(X^{m+n} - p X^{m}) =  p X^{2m-n} \in I^2.$$ Thus as $ 2m-n \leq 2n, pX^{2n} \in I^2 $. Cosequently $$ p( X^{2n} - 2p X^{n} +p^2)- p X^{2n}  = p^3 - 2 p^2 X^{n} \in I^2.$$ 
                               Now note that in Step 1, we have proved that $p X^{m-n}\in I$. Hence $ p^2 X^{n} \in I^2$ since $ 2(m-n) \leq n$. Consequently 
                                $$ (p^3 - 2p^2 X^n) + 2 p^2 X^n = p^3 \in I^2.$$    Thus Step 2  is proved.\\                             
                               By Step 1 and Step 2, it is immediate that $ I^2 \cap R \neq (I\cap R)^2 $ i.e., $I$ is not power stable. Further, note that for any $ n \geqq 2$, $ p^{2n+1} \in I^{n+1}$, hence $ (I\cap R)^n / (I\cap R)^{n+1} \longrightarrow I^n/I^{n+1} $ is not a monomorphism. Thus the ideal $I$ is not ultimately power stable.\\
                              
                               Remark 3.17. In the above example, radical of $I$ is $ id.(X, p)$, a maximal ideal in $R[X]$. Thus $I$ is a primary ideal. Hence primary ideals need not be power stable.\\
                               
                               Theorem 3.18. Let $R$ be an integral domain. If every ideal in $R[X]$ is power stable, then $R$ is a field.\\
                               
                               Proof. Assume $R$ is not a field. Let $ \wp$ be a prime ideal in $R$. By our assumption, every ideal in $ R_{\wp}[X]$ is power stable. Thus, it is clear that to prove the result we can assume that $R$ has a unique maximal ideal. Let $ \textbf{m}$ be the maximal ideal in $R$. By assumption, for any  $ \lambda \in \textbf{m} - \textbf{m}^2,  J = id.(X^2 - \lambda , \lambda X) $ is power stable. Put $ I = J \cap R$. If $ a\in I$, then      \begin{eqnarray*} a & = & f(X) (X^2 - \lambda) + g(X) \lambda X \hspace { 1 cm } ( f(X), g(X) \in R[X])\\
                               	\Rightarrow a & = & - f(0) \lambda \text{  putting  X= 0  } \\  \Rightarrow I & \subset & R \lambda \\ \Rightarrow I & = & I_1 \lambda 
                               \end{eqnarray*}   where  $ I_1 = \{ b\in R : b\lambda \in I\} $. We, now, consider two cases. \\
                               
                               Case 1.  $ I_1 \subset \textbf{m}$.\\
                               
                               Note that  $ \lambda^3 = \lambda (X^2 - \lambda)^2  + (\lambda X)^2 - \lambda X^2 (X^2 - \lambda).$     Hence  \begin{eqnarray*} \lambda^3 \in J^2 \cap R & = & I_1^2 \lambda^2  \\  \Rightarrow  \lambda^3 & = & b \lambda^2 \hspace { .5 cm} ( b \in I_1^2)\\  \Rightarrow \lambda & = & b \in I_1^2 \subset \textbf{m}^2. \end{eqnarray*}
                                                    
                                 This, however, is not true by choice of $ \lambda$. Therefore Case 1 does not occur.\\
                                 
                                 Case 2. $ I_1 = R$. \\
                                 
                                 In this case $J^2 \cap R = R\lambda^2$. Thus, as $\lambda^2 \in J^2$, we have $$ \lambda^2 = (X^2- \lambda)^2 a(X) + (\lambda X)^2 b(X) +\lambda X(X^2 - \lambda)c(X)$$ \\where $ a(X), b(X), c(X) \in R[X]$. Thus $a(0) = 1$. Hence $$ 0 = X^4.  \overline{a(X)} \text{ in } R/(\lambda)[X]$$   This, however, is not true as $ a(0) = 1$. Consequently $R$ is a field.\\
                                 
                                 Remark 3.19. The above Theorem is not valid even if we replace $R[X]$ by a faithful flat extension ring of $R$ e.g., it can be easily checked that for $ \mathbb{Z} \subset \mathbb{Z}[i]$, every ideal of $ 
                                 \mathbb{Z}[i]$ is power stable in $\mathbb{Z}$.\\
                                 
                                 Theorem 3.20. Given any power stable ideal $I$ in $R[X]$, $I\cap R = J$,there exists a power stable ideal $I_1 \subset R[X]$ maximal with respect to the property $ I \subset I_1, I_1 \cap R = J $. Further, if $J$ is prime, then whenever $ f(X). g(X) \in I_1, f(X) \notin I_1, g(X) \notin I_1$, then  $ (I_1 + id.(f(X)) \cap R = J$ or  $ (I_1 + id.(g(X)) \cap R = J$.\\
                                 
                                 Pf. Let $$ \mathfrak{S} = \{ \mathfrak{a} \subset R[X], \text{ power stable ideal } : I \subset \mathfrak{a}, \mathfrak{a} \cap R = J  \} .$$ The set $ \mathfrak{S}$ is partially ordered with respect to containment. Let $ \{\mathfrak{a}_i, i \in I \}$ be a chain in $ \mathfrak{S}$. Put $ \mathfrak{b} = \cup \mathfrak{a}_i$, then $ \mathfrak{b}$ is an ideal in $R[X]$, and $ \mathfrak{b}\cap R = \cup (\mathfrak{a}_i \cap R) = J.$ Moreover  $$ \mathfrak{b}^n \cap R = \cup (\mathfrak{a}_i^n \cap R) = J^n.$$  Hence $\mathfrak{b}$ is an upper bound of the chain $ \{\mathfrak{a}_i, i \in I \}$  in $ \mathfrak{S}$. Therefore by Zorn's Lemma, $ \mathfrak{S}$ has a maximal element $ I_1$ (say). We shall now show that if $J$ is prime then the last part of the statement holds. Assume the polynomials $ f(X) \notin I_1, g(X) \notin I_1$, but $ f(X).g(X) \in I_1$. If  $ J \subsetneqq ( I_1 + id.(f(X)) \cap R $ and  $ J \subsetneqq ( I_1 + id.(g(X)) \cap R $, then there exist $ a(X), b(X) \in I_1$ such that $ a(X) + f(X) h(X) \in (R- J), b(X) + g(X) h_1(X) \in (R-J)$. Then $ (a(X) + f(X) h(X)) (b(X) + g(X) h_1(X)) \in I_1 \cap R = J$. This cotradicts the fact that $J$ is prime. Hence the result holds.\\                                 
                                 
                                 Theorem 3.21. (i) Let $ A \subset B$ be rings, and  $ I $, an ideal in $B$, which is an ultimately power stable in $A$. Then if $ I \cap A$ is finitely generated projective $A-$module of rank one, $I$ is power stable.\\
                                 
                                 (ii)  Let $ A \subset B$ be Noetherian rings, and  $ I $ is an ideal in $B$ which is ultimately power stable in $A$. Then if $ I \cap A = J$ is a radical ideal which is generated by a regular $A-$ sequence, it is power stable.\\
                                 
                                 Proof. (i) Clearly, for all $ n \geq 1$, we have $ (I\cap A)^n \subset I^n \cap A$. Thus if $I$ is not power stable there exists $m > 1$ largest 
                                 	
                                 such that \begin{eqnarray*}(I\cap A)^m & \subsetneqq & I^m \cap A \\ \Rightarrow    (I\cap A)^{m+1} &  \subset &  (I\cap A) (I^m \cap A)\\ &  \subset  & I^{m+1}  \cap A  \\                                  	  \end{eqnarray*}
                                 Note that $ I^{m+1} \cap A = (I\cap A)^{m+1}$ by maximality of $m $. Hence 
                                   \begin{eqnarray*} (I\cap A )^{m+1} & = & (I \cap A) ( I^m \cap A) \\ \Rightarrow (I \cap R)_{\wp}^m  & = & (I^m \cap A)_{\wp} \hspace{.5 cm}  \forall  \wp \in Spec(A) \text { since }  (I\cap A)_{\wp}  \text{ is free of rank 1 }  \\ \Rightarrow  (I \cap A)^m & = & I^m \cap A.                              
                                 \end{eqnarray*} This is in contradiction to the choice of $ m $. Hence $I$ is power stable in A.\\
                                 
                                 (ii) Let $ \{g_1,\cdots, g_n \} $ be a regular sequence in $A$  generating the ideal $J$. By  [5, Theorem 2.1], there exists an isomorphism from the graded polynomial ring  $(A/J)[X_1, \cdots, X_n], $  where $ X_1, \cdots, X_n$ are indeterminates over $A/J$, to the graded ring $ G_J(A)$ which maps each $X_i, i=1,2, \cdots , n$ to the image of $g_i$ in $G_J(A). $ Note that, as $J$ is a radical ideal, the ring $(A/J)[X_1, \cdots, X_n] $ is  without nilpotents. Thus $G_J(A)$ does not have nilpotents. Consequently  $I$ is power stable by  Theorem 2.12 (ii). \\  
                                 
                                 We want to further understand the relationship in ultimate power stability and power stability of an ideal. Our next theorem in this regard uses a result in [4] which we record for convenience :[4,Theorem 2, pp. 156] Let $A$ be a Noetherian ring, and $ J \subset I$ be ideals in $A$ where $I$ has a non-zero divisor. Then $ J$ is a reduction of $I$ if and only if $ IK =JK$ for an ideal $K$, which contains a non-zero divisor. \\        
                                                                  
                                Theorem 3.22.   Let $ A \subset B$ be  integral domains where $A$ is Noetherian local ring. Assume
                                 $ I $ is an ideal in $B$, which is an ultimately power stable in $A$ and $( I^l \cap A) = (I\cap A)^l$ for all $ l > m$,but $ (I^m \cap A) \neq (I\cap A)^m$. Then for $ 1 \leq t \leq m$, we have\\
                                
                                (i)   $ \forall  1 \leq t \leq (m+1-t), (I^{m+1-t}\cap A)^r = ((I\cap A)^{m+1-t})^r$ for some $ r \geq 1$.\\
                                                    
                                (ii) $ \forall   t > (m-t+1)$, we have $ (I\cap A)^{m-t+1} (I^{m-t+1}\cap A)^r = ((I\cap A)^{m-t+1})^{r+1}$ for some $ r \geq 1$\\
                                
                               Proof. We shall prove the result in steps:\\
                               
                               Step 1. $ (I\cap A)^t (I^{m+1-t} \cap A) = (I \cap A)^{m+1} $ for all $ 1 \leq t \leq m$.\\
                               
                              From the proof of Theorem 3.21(i), it follows that the assertion is true for $t = 1.$ Hence assume that for $ 1 \leq s < m$, we have $$ (I\cap A)^s ( I^{m+1-s}\cap A)= (I\cap A)^{m+1}$$ Then as $$ (I\cap A)^{s+1} ( I^{m-s}\cap A)= (I\cap A) (I\cap A)^s ( I^{m-s}\cap A) \subset (I\cap A) (I^m\cap A),$$ we get $$ (I\cap A)^{s+1} ( I^{m-s}\cap A) \subset (I\cap A)^{m+1} = (I\cap A)^{s+1} ( I \cap A)^{m-s}.$$ Now, as $ (I\cap A)^{m-s}\subset (I^{m-s} \cap A)$, we conclude $$ (I\cap A)^{s+1} (I^{m-s}\cap A) = (I\cap A)^{m+1}.$$ Cosequently Step 1 is proved.\\
                              
                              Step 2. $ \forall  1\leq t \leq (m+1-t)$, we have $  (I^{m+1-t}\cap A)^r = ((I\cap A)^{m+1-t})^r$ for some $ r \geq 1$\\
                              By Step 1, we have $$ (I\cap A)^t (I^{m-t+1}\cap A) = (I\cap A)^{m+1}.$$ Thus $$ (I\cap A)^t (I^{m-t+1}\cap A) = ( I\cap A)^t (I\cap A)^{m-t+1}.$$ Therefore by  [4,  Theorem 2,pp.156 ], we conclude that  $(I\cap A)^{m-t+1} $ is a reduction of $(I^{m-t+1}\cap A)$. Hence there exists $ l \geq 1$ such that  $$(I\cap A)^{m-t+1} (I^{m-t+1}\cap A)^l =(I^{m-t+1}\cap A)^{l+1}.\ldots  (1)$$ As there are only finite number of $t's$ in question, we can choose common l, although we do not need this. As $ m+1-2t \geq 0$, multiplying the equation in Step 1,  by $ (I\cap A)^{m+1-2t}$, we get  $$ (I\cap A)^{m-t+1}(I^{m-t+1}\cap A)= (I\cap A)^{2(m-t+1)}  \ldots (2)$$ From the equation (2), for any $ l \geq 1$, we have $$ (I\cap A)^{m-t+1}(I^{m-t+1}\cap A)^l= (I\cap A)^{(l+1)(m-t+1)} . $$  Hence using the equation (1), we get $$ (I\cap A)^{(m-t+1)(l+1)} = (I^{m-t+1}\cap A)^{l+1}.  $$ Thus Step 2 is proved.\\
                              
                              Step 3.  $ \forall   t > (m-t+1)$, we have $ (I\cap A)^{m-t+1} (I^{m-t+1}\cap A)^r = ((I\cap A)^{m-t+1})^{r+1}$ for some $ r \geq 1$\\ 
                              
                               By Step 1, we have $$ (I\cap A)^{2t} (I^{m-t+1}\cap A) = (I\cap A)^{m+t+1} =(I\cap A)^{2t} (I\cap A)^{m-t+1} .$$\\ Therefore by  [4,  Theorem 2,pp.156 ], we conclude that  $(I\cap A)^{m-t+1} $ is a reduction of $(I^{m-t+1}\cap A)$. Hence $$ (I\cap A)^{m-t+1} (I^{m-t+1}\cap A)^r = (I^{m-t+1}\cap A)^{r+1}$$ for some $ r \geq 1.$ Hence the theorem.\\  
                               
                               Remark 3.23. In the proof of the above theorem, the only fact we have used is that $( I \cap A) (I^m \cap A ) = (I \cap A)^{m+1}$.\\
                              
                              We, now, note a general observation which strengthens Theorem 3.22(i).\\
                              
                              Theorem 3.24.  If $ P \subset K$ are ideals in a ring $A$ such that $ P^l = K^l$ for some $ l \geqq 1$, then $ P^t = K^t$ for all $ t \geqq l$.\\
                              
                              Proof. Let $ r \leqq l$. Then $$ P^r K^{l-r}  \subset K^l  = P^l = P^r P^{l-r} \subset P^r K^{l-r}.$$ Hence $ P^r K^{l-r} = P^l =K^l$ for all  $ r \leqq l$. Therefore $$ P^{l+1} \subset K^{l+1} = K P^l = K P^{l-1} P = P^l P= P^{l+1}.$$ Hence $ K^{l+1} = P^{l+1}$. Consequently, by by induction  it follows that $ P^ t = K^t $ for all $ t \geqq l$.\\
                              
                              Remark 3.25. In the above theorem, if $A$ is a Noetherian integral domain, then  by [4, Theorem 2, pp. 156 ], under the given conditions, $K$ is a reduction of $P$.\\

                                It would be interesting to know the answer to the following : \\
                            
                          {\bf{Question}}1. Let $R$ be a Noetherian domain of dimension 1. Does there exists a characterization of power stable ideals in $R[X]$.  \\
                          
                         {\bf{Question}} 2. For an integral domain $R$ of dimension n, does there exist $ m \geq 2$ such that for any ideal $I \subset R[X], I^t \cap R = (I\cap R)^t$ for all $ 1 \leq t \leq m$ implies $ I^t \cap R = (I\cap R)^t$ for all $ t \geq 1$.\\

                             \begin{center}  ACKNOWLEDGEMENGT
                             	
                             \end{center}
                             I am thankful to Stephen McAdam and Mosche Roitman for some useful e-mail exchanges. I am also thankful to Melvin Hochster for the examples and for his very prompt responses to my queries. \\
                             
                              \begin{center}  REFERENCES                              	
                              \end{center}
                              
                              1. M.F.Atiyah , I.G. Macdonald, Introduction to Commutative Algebra, Addison- Wesley Publ. Co., 1969.\\
                              
                              2. S.M.Bhatwadekar, Pramod K. Sharma, Unique Factorization and Birth of Almost Primes, Communications in Algebra, 33:43-49, 2005.\\           
                                                          
                              3. Irving Kaplansky, Commutative Rings, The University of Chcago Press, 1974.\\
                              
                              4. D.G.Northcott and D.Rees, Reductions of Ideals in Local Rings, Proc. of Camb. Phil. Soc., 1954 (50), 145-158  \\
                              
                              5. David Rees , The Grade of an Ideal or Module, Proc. of Camb. Phil. Soc., 1957 (53), 28-42.\\
                              
                             6. Pramod K. Sharma, On Power Stable Ideals, arXiv:0705,1286v1[math.Ac] 9th May, 2007.

                     \end{document}